\documentclass[conference]{IEEEtran}
\IEEEoverridecommandlockouts
\usepackage{cite}
\usepackage{amsmath,amssymb,amsfonts, dsfont}
\usepackage{algorithm, algpseudocode, booktabs}
\usepackage{graphicx, subcaption}
\usepackage{float}
\usepackage{textcomp}
\usepackage{xcolor}
\def\BibTeX{{\rm B\kern-.05em{\sc i\kern-.025em b}\kern-.08em
    T\kern-.1667em\lower.7ex\hbox{E}\kern-.125emX}}
\begin{document}
\title{Implementing Recycling Methods for Linear Systems in Python with an Application to Multiple Objective Optimization
}

\author{
\IEEEauthorblockN{Ainara Garcia}
\IEEEauthorblockA{\textit{College of Engineering, Computing \& Applied Sciences} \\
\textit{Clemson University}\\
Clemson, SC, USA \\
ainarag@clemson.edu}
\and
\IEEEauthorblockN{Sihong Xie}
\IEEEauthorblockA{\textit{Thrust of AI} \\
\textit{HKUST-GZ}\\
Guangzhou, China \\
xiesihong1@gmail.com}
\and
\IEEEauthorblockN{Arielle Carr}
\IEEEauthorblockA{\textit{Computer Science \& Engineering} \\
\textit{Lehigh University}\\
Bethlehem, PA, USA \\
arg318@lehigh.edu}
}

\maketitle
\begin{abstract}
Sequences of linear systems arise in the predictor-corrector method when computing the Pareto front for multi-objective optimization. Rather than discarding information generated when solving one system, it may be advantageous to recycle information for subsequent systems. To accomplish this, we seek to reduce the overall cost of computation when solving linear systems using common recycling methods. In this work, we assessed the performance of recycling minimum residual (RMINRES) method along with a map between coefficient matrices. For these methods to be fully integrated into the software used in Enouen et al. (2022), there must be working version of each in both Python and PyTorch. Herein, we discuss the challenges we encountered and solutions undertaken (and some ongoing) when computing efficient Python implementations of these recycling strategies. The goal of this project was to implement RMINRES in Python and PyTorch and add it to the established Pareto front code to reduce computational cost. Additionally, we wanted to implement the sparse approximate maps code in
Python and PyTorch, so that it can be parallelized in future work.

\end{abstract}

\section{Introduction}

Issues of fairness often arise in graphical neural networks used for misinformation detection. Accurately reducing unfairness when satisfying more than one constraint can come at a significant computational cost. Multi-objective optimization (MOO) methods allow us to explore different trade-offs by generating a set of solutions, known as the Pareto front \cite{b1}. Traditional first-order MOO methods, such as multi-gradient descent (MGD)~\cite{b2}, are computationally expensive on large graphs due to the gradient computation at each iteration, and may require many iterations to generate a Pareto front. Using the predictor-corrector method introduced in \cite{b3}, and iterative methods like MINRES\cite{b4} and Conjugate Gradient (CG)\cite{b5}, the linear systems arising in the predictor step can be solved efficiently. 
 
We specifically focus on the iterative method MINRES\cite{b4} and recycling MINRES (RMINRES)\cite{b6} for solving the linear systems arising in the predictor step. In \cite{b3}, the authors found that the computational cost associated with solving this sequence of systems can be reduced by dramatically limiting the maximum number of iterations when using MINRES, while still achieving an accurate Pareto front \cite{b3}. In the current work, we are interested in assessing the performance of recycling strategies across the sequence of linear systems to further reduce the total number of iterations.  In this project, we also consider a mapping strategy between the coefficient matrices that may allow us to reuse previous solutions. Specifically, we explore the potential for a parallel version of the preconditioner update reported in \cite{b7} to further reduce computational costs associated with many preconditioners for long sequences of systems arising in the predictor step. In particular, we seek to implement both RMINRES and the preconditioner updates in Python and PyTorch\cite{b8}, using the respective language and libraries when computing the Pareto front's predictor-corrector methods. We provide preliminary results related to our work, noting that the implementation is still in progress. Aiming to help other groups working on the same topic, we also detail the complications encountered.

\section{Preliminaries}
\label{sec:prelim}

\subsection{Predictor Corrector Method}

\begin{figure}[h!]
    \centering
    \includegraphics[scale=0.25]{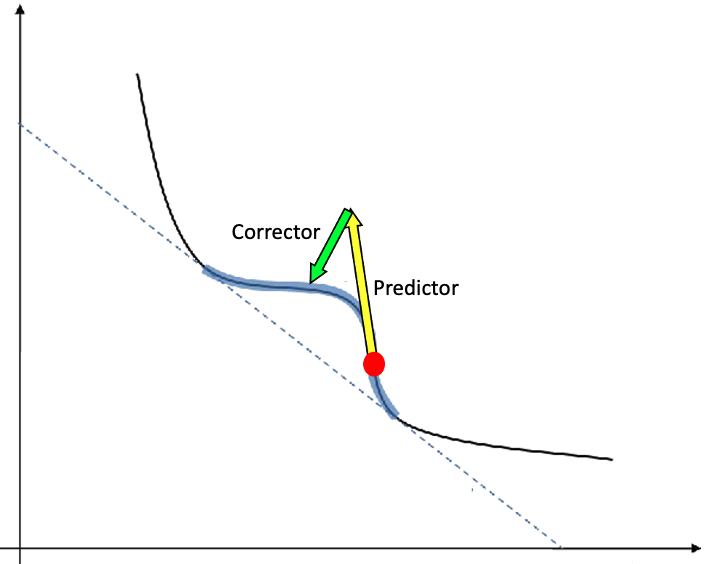}
    \caption{The predictor step (yellow) and corrector step (green) given an initial point (red).}
    \label{fig:PC}
\end{figure}

Given a point of a function, the predictor-corrector method allows us to approximate the value of that function at a nearby point. This method consists of two steps (see Figure \ref{fig:PC}). First, in the predictor step, we determine an approximate direction to a neighboring point, and then move along that direction based on a predetermined, fixed step size. Then, we refine this initial approximation in the corrector step. 

To determine the direction in which to move in the predictor step, we solve for the vector $v$ in Equation 1,
\begin{equation}
H(x_0^*)v = \nabla f(x_0^*) \beta
\label{eq:H}
\end{equation}
where $x_0^*$ is a starting Pareto optimal point, and $H(x_0^*)$ and $\nabla f(x_0^*)$ are the Hessian and Jacobian of our loss functions at $x_0^*$, respectively; and $\beta$ is a weighting vector. Intuitively, this approach works by finding a suitable direction to move along the tangent plane of $f$ at $x_0^*$ as depicted in Figure \ref{fig:PC}. While more information on this method can be found in \cite{b3}; this report focuses on the reducing costs associated with solving the sequence of systems of the form (\ref{eq:H}) that arise in the predictor step.
Directly solving the linear systems of the form (\ref{eq:H}) requires $O(n^3)$ time, which is generally expensive, even for modest matrix dimensions, $n$. Furthermore, in this application, $H$ is not explicitly known, so a matrix-free method, such as MINRES is required\cite{b3}. Thus, we consider solvers whereby the solution is approximated iteratively (i.e., iterative solvers). Since the coefficient matrix $H$ is symmetric, we use a symmetric solver such as MINRES, which requires only a linear operator to compute a matrix-vector (the main computational cost in such iterative methods). Instead of starting the computation of MINRES from scratch for each new linear system, the recycling approach allows us to reuse information from previous solutions to accelerate convergence and reduce computational costs by reducing the total number of iterations to solve the sequence of linear systems. To accomplish this, we consider the RMINRES method\cite{b6}.

\section{Iterative Methods Solving Linear Systems}
Iterative methods compute an approximate solution to a linear system by repeating a series of calculations to update the initial guess until convergence is achieved. Unlike direct methods, which compute an exact solution, iterative methods progressively refine an approximate solution until the accuracy is below the pre-determined tolerance. The processes includes taking an initial guess, which will be the starting vector for the rest of the iteration process. 

The MINRES algorithm \cite{b4} is a Krylov iterative method used to solve systems of linear equations when the coefficient matrix is symmetric and possibly indefinite. MINRES was developed as an alternative to other well-known symmetric solvers, such as the CG algorithm\cite{b5}, which requires the coefficient matrix to be symmetric, positive, and definite. A distinguishing features of MINRES is its ability to handle indefinite matrices. In contrast to positive definite matrices, which have all positive eigenvalues, indefinite matrices might have both positive and negative eigenvalues. 

Central to MINRES is the emphasis on minimizing the residual vector norm, $\|r_k\|_2$, at each iteration, where the residual of an approximate solution $x_k$ is defined as $r_k = Ax_k - b$. The Krylov subspace of dimension $k$ is defined as $\mathcal{K}^k(A, r_0) = {span}\{r_0, Ar_0,\,...,\,A^{k - 1}r_0\}$, where $A$ is our coefficient matrix from the linear system $Ax =b$, $r_0 = b - Ax_0$ is the initial residual, and $x_0$ is the initial guess. $\mathcal{K}^k(A, r_0)$ is expanded by one vector with each iteration of MINRES. The $k^{th}$ iteration $k$ of MINRES computes $z_k \in \mathcal{K}^k(A, r_0)$ and updates $x_k$ as $x_k = x_0 + z_k$ such that the residual norm $\|r_k\|_2 = ||Ax_k - b\|_2$ is minimized. This continues until $\|r_k\|_2$ is within a predefined tolerance. More information on MINRES and other Krylov methods can be found in \cite{b12}. This process ensures that convergence can be achieved in $n$ iterations, but ideally an accurate-enough solution can be achieved in $m \ll n$ iterations. To further reduce the number of iterations even more, we can recycle some of the Krylov vectors generated for a previous system to start the next call to the solver for a subsequent system rather
than starting with just one vector. This new set of
vectors is stored in the columns of a matrix we refer to as $U$, called the recycle space \cite{b6}. 

\begin{figure}
    \centering
    \includegraphics[width=9cm]{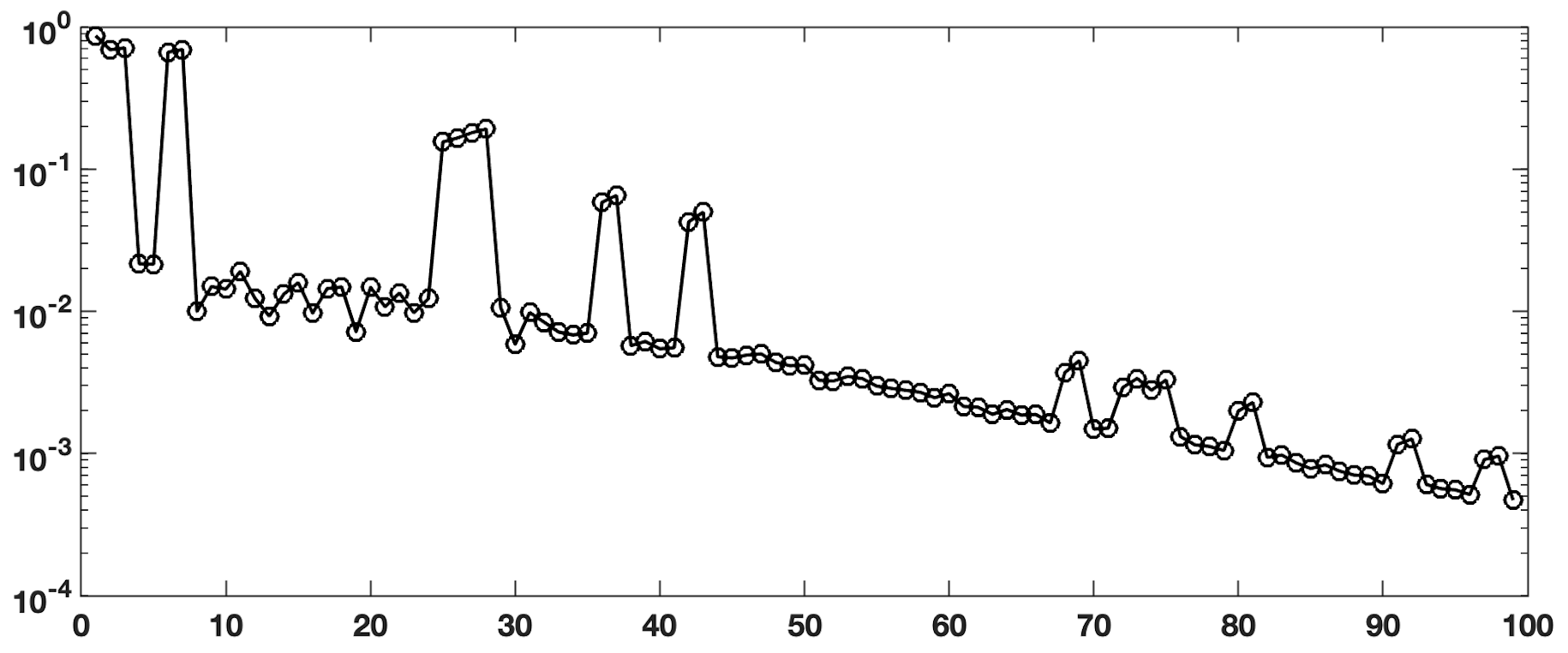}
    \caption{Normwise difference between each pair of matrices in a sequence of 100 systems extracted from \cite{b1}. Note that the y-axis shows $log\left(\|A_k-A_{k-1}\|_F\right)$.}
    \label{fig:normDiff}
\end{figure}

For linear systems that are close to one another, recycling subspaces (i.e., $U$) can be advantageous to reduce the total number of iterations of the iterative solver across the entire sequence.  We provide iteration comparisons in the next section, but refer to Figure \ref{fig:normDiff} to demonstrate how closely related the systems are using 100 matrices extracted from \cite{b3}. We plot the norm-wise difference between each pair of matrices $\|A_k - A_{k-1}\|_2$, for $k = 2, 3, \dots, 100$. The dimension of these systems is $n = 1234$. The peaks generally occur at each new iteration of the predictor-corrector method; as noted in \cite{b3} and references within. Even in cases of these ``larger" peaks, the differences between the matrices are quite small, motivating our conjecture that the subspaces generated by one call to MINRES will be a close approximation to the subspace for the next call (and thus will serve as a good recycle space for RMINRES).

\section{Recycling Strategies}
The iterative nature of the MINRES algorithm enables it to handle large-scale linear systems that might be computationally expensive to solve using direct methods. By gradually refining the approximate solution with each iteration, MINRES offers an attractive compromise between accuracy and computational efficiency. Many algorithms also use preconditioners as part their computation to further accelerate convergence when the matrix is ill-conditioned, that is $\|A\|\|A^{-1}\| \gg 1$. However, because of the computational cost associated with computing a good preconditioner, we can instead recycle preconditioners by updating a previous preconditioner to then be reuse for the next system.  


Preconditioners are used in iterative methods to speed up the convergence of linear systems when the coefficient matrix is ill-conditioned. By transforming the original system into a better-conditioned one, preconditioners reduce the number of iterations required to reach an accurate solution. More details on preconditioning techniques can be found in \cite{b9} and references therein. Preconditioners can enhance the robustness of iterative solvers, making them more tolerant to ill-conditioned or difficult-to-solve problems. They can handle cases where the original system matrix is poorly conditioned or close to singular. However, they can be expensive to compute and so we consider preconditioner updates, specifically the sparse approximate map\cite{b7}. 

\begin{figure}[h]
    \centering
    \includegraphics[width=9cm]{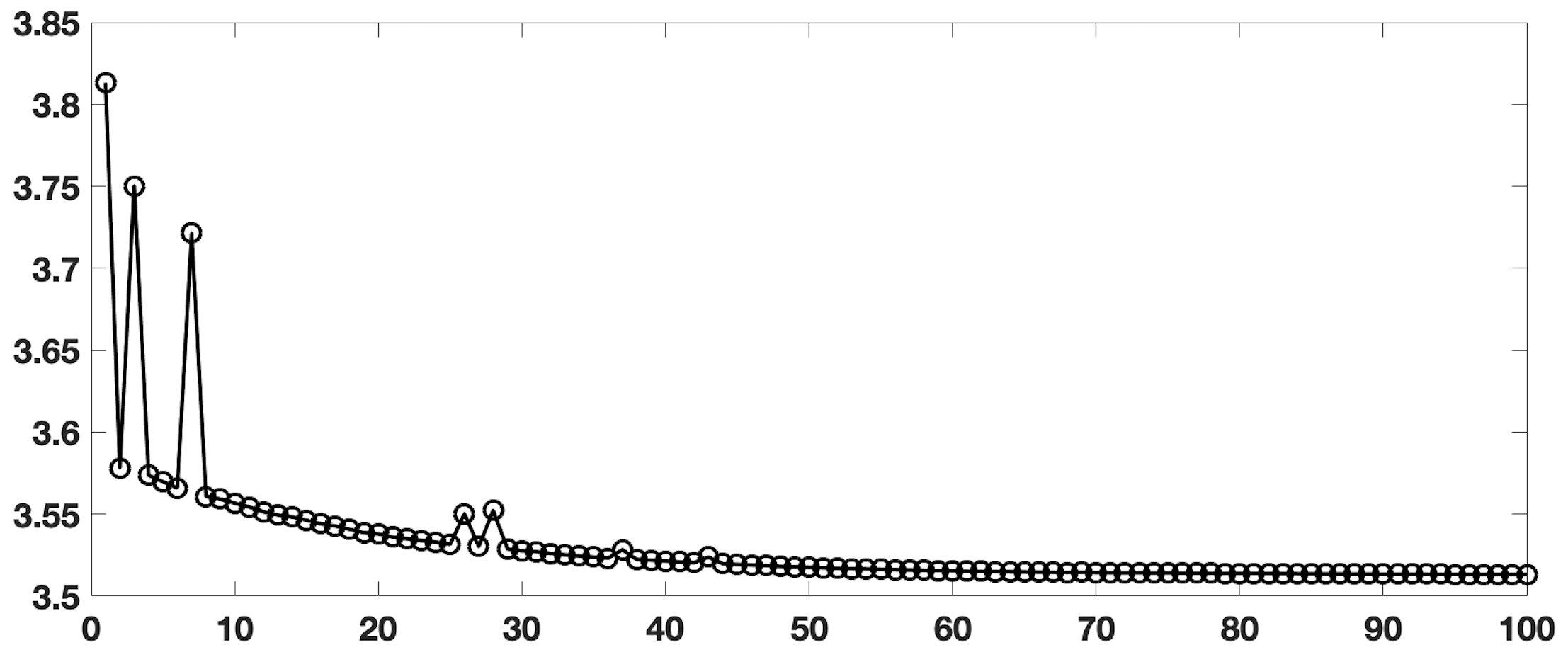}
    \caption{Condition number (y-axis) of matrices from a sequence of 100 systems extracted from \cite{b1}.}
    \label{fig:cond}
\end{figure}

The authors in \cite{b7} present the central concept behind SAMs: the mapping of one matrix to another for which we have a good preconditioner (i.e., one that facilitates fast convergence of the iterative solver) \cite{b7}. For instance, consider a matrix $A_k$ in a sequence of matrices, we want to compute a map $N_k$ to another matrix $A_0$ with the constraint of an imposed sparsity pattern to limit the cost of computing and applying $N_k$. In particular, the equation $A_kN_k = A_0$ expresses the exact mapping, where $A_k$ represents the current system matrix, $A_0$ is the target system matrix, and $N_k$ is the map. The objective is to find an approximate $N_k$ such that it approximates $A_0$ well enough. Due to the iterative nature of performing SAMs, we can reduce the computational cost of solving linear systems and increase efficiency of our implementation by exploring the potential of a parallel implementation.  As shown in Figure \ref{fig:cond}, the condition numbers (i.e., $\|A\|\|A^{-1}\|$) of our systems is relatively low, so while the use of a preconditioner in this application was ultimately not necessary, we still consider a Python implementation for future use in cases where the matrices are, in fact, ill-conditioned. In section \ref{sec:SAMS}, we give the explicit mathematical expressions that define a SAM and discuss the corresponding implementation challenges and (proposed) solutions. 

There is often a trade-off between accuracy and the cost of (preconditioned) iterative methods. Iterative methods typically allow control over the accuracy by specifying a convergence criterion (e.g., a tolerance level or maximum iterations allowed). Achieving higher accuracy may require more iterations, leading to increased computational cost. Choosing an appropriate iterative method, setting suitable convergence criteria, and considering the resources available are crucial in striking the right balance between accuracy and cost.
Krylov methods, such as MINRES and RMINRES, are particularly useful in our case since these methods do not require the matrix to be explicitly defined, and are often referred to as `matrix-free methods'\cite{b10}. In our case, $H$ in \ref{eq:linear_system} is not known explicitly, only its action on a vector.

\begin{figure}
    \centering
    \includegraphics[width=9cm]{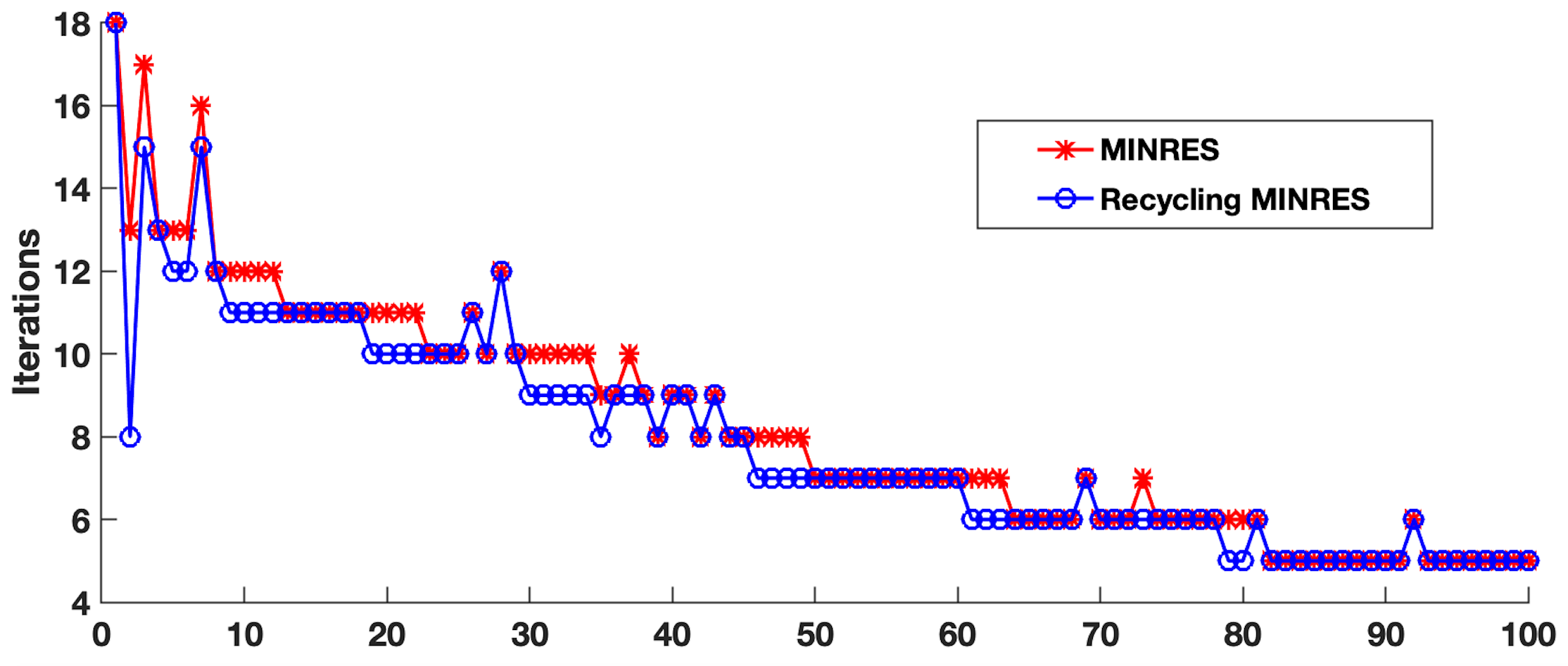}
    \caption{Iterations for MINRES and RMINES for a sequence of 100 systems extracted from \cite{b1}.}  \label{fig:MINRESvsRMINRES}
\end{figure}

The goal of this project is to implement RMINRES in Python and add it to the established Pareto front code used in \cite{b3} to reduce the computational cost. Additionally, the project aimed to implement the sparse approximate maps code in Python so that in future work, it can be parallelized. We ran initial experiments using a MATLAB implementation of MINRES and RMINRES on the same sequence of linear systems as in Figures \ref{fig:normDiff} and \ref{fig:cond}. In Figure \ref{fig:MINRESvsRMINRES}, we compare the iterations for the two methods.  
We see that RMINRES reduces total iterations from 828 to 796, and we observe the largest reduction in iterations for the linear systems arising earlier in the sequence.
These results motivate our implementation of RMINRES in Python/PyTorch as the cost of the Hessian-vector products impose a significant computational computational bottleneck.  As the size of the systems grow (i.e., as the number of constraints increase in the optimization), such reductions become critical.

\section{Implementation}
It is important to have iterative solver methods in an industry-standard computer science language. For the code implemented in this project, the Python language and several libraries were used. NumPy, KryPy, and SciPy\cite{b11} packages played a significant role in these implementations. PyTorch is an integral tool in machine learning and parallel computing. It allows for scalable  distributed training and performance optimization in research. However, like all tools, it has limitations. Initially, when implementing RMINRES in PyTorch, an error within the implementation occurred. However we discovered that the sparse tensor data structure and operations required for this project, are not developed enough for the scale of our research. The PyTorch linear algebra solver, torch.linalg.solve, often used to find the solution of systems of linear equations, is compatible with float, double, cfloat, and cdouble dtypes \cite{b8}. This limitation in the PyTorch operations caused errors in the implementation when solving the incomplete Cholesky factorization of coefficient matrix. This implementation in PyTorch presented many challenges when trying to compute matrix-vector multiplication. Therefore, the alternative choice was to implement the code in Python initially, then integrate it into PyTorch from Python. 
Python has many packages capable of handling large-scale computing, matrices, linear algebra, and iterative methods. 

KryPy, a Krylov subspace methods package \cite{b11}, includes many of the functions called into our implementations. Among these packages is krypy.recycling, which allows for the computation of solutions of linear systems using recycling. Calling krypy.recycling allows us to build the $U$ subspace and return it to be implemented in the next call to RMINRES. Many combinations of KryPy functions were attempted before coming to the conclusion that KryPy's recycling functions are limited to the number of iterations allowed, and many times does not converge within those iterations for the computations done in this project. The SciPy library has a package that allows for the computation of sparse linear systems, for example the compressed sparse column format used in this project, as well as functions used to process matrices, like returning a sparse matrices from diagonals and performing linear algebra operations. NumPy is commonly used for scientific computation in Python. While these libraries were important in the development of the code for recycling strategies, their data structures are sometimes not compatible with PyTorch, and therefore require significant modification. 


\section{Workflow and Steps}
As the goal for the project is to integrate these commonly used methods to improve performance of the predictor-corrector method, we began with sequential implementations of our methods. When referencing MINRES, it is important to again emphasize that it begins with an initial vector, while RMINRES begins with many vectors, or a subspace generally defined by vectors from the Krylov space built for the previous linear system in the sequence. This often leads to a reduction in solver iterations compared with standard MINRES (as shown in Figure \ref{fig:MINRESvsRMINRES}). The reader can refer to \cite{b10} for the detailed algorithm for RMINRES. 

 \begin{algorithm}
 \small
 \caption{Conjugate Residual Method for $H\textbf{v} = \textbf{b}$}\label{alg:minres}
 \begin{algorithmic}[1]
 \State $\textbf{r}^{(0)} = \textbf{b}-H\textbf{v}^{(0)}$.
 \State $\textbf{p}^{(0)} = \textbf{r}^{(0)}$.
 \State $i = 0$.
 \While{$i < maxIter$ and $\|\textbf{r}^{(i)}\|> tol$}
 \State $\alpha_{i} = \frac{(H\textbf{r}^{(i)})^\top \textbf{r}^{(i)}}{\|H\textbf{p}^{(i)}\|^2} $.
 \State $\textbf{x}^{(i+1)} = \textbf{x}^{(i)} + \alpha_{i}\textbf{p}^{(i)}$.
 \State $\textbf{r}^{(i+1)} = \textbf{r}^{(i)} - \alpha_{i}H\textbf{p}^{(i)}$.
 \State $\beta_i = \frac{(H\textbf{r}^{(i+1)})^\top\textbf{r}^{(i+1)}}{(H\textbf{r}^{(i)})^\top\textbf{r}^{(i)}}$.
 \State $\textbf{p}^{(i+1)} = \textbf{r}^{(i+1)} + \beta_{i}\textbf{p}^{(i)}$.
 \State $i = i + 1$.
 \EndWhile
 \end{algorithmic}
 \end{algorithm}

\subsection{Recycling Minimum Residual (RMINRES)}
In our case, several functions are required to implement RMINRES. We also need to calculate the $k$ eigenvalue-eigenvector pairs, referred to as Ritz pairs (see e.g., \cite{b8}), that correspond to the smallest eigenvalues of the coefficient matrix. Then, we generate a sparse matrix $A$ and right-hand-side $b$ from test matrices. Algorithm 1 gives the conjugate residual method \cite{b10}, an equivalent version of MINRES; we use this algorithm for simplicity, but refer the reader to \cite{b4} for an in depth discussion of the MINRES algorithm.

RMINRES function requires several initial inputs: $A$ (from $Ax=b$), $b$ (again from $Ax=b$), an initial guess ($x_0$), a convergence tolerance, the maximum number of iterations, a preconditioner (if available/desired), and the recycle space $U$. The code checks if a preconditioner and a recycle space were provided, then computes the initial residual (corresponding to line 1 of Algorithm 1), $r = b - Ax_0$. It then computes the norm of the initial residual and checks if the initial guess $x_0$ is a converged solution. If yes, it returns the result (but generally this is not the case). If a recycle space $U$ is provided, it orthogonalizes it against the current approximation to the solution (current approximation computed on line 6 of Algorithm 1), updates the initial guess and residual, and checks for convergence again. The symmetric tridiagonal matrix $T$ generated by the Lanczos method\cite{b12} is not given in Algorithm 1, but we note that this is the procedure used by MINRES and RMINRES to build an orthogonal basis for the Krylov space using successive matrix-vector products.  These are fairly easy to implement in Python, however we must compute the QR factorization of $T$ to compute the update to the approximate solution. 

The code continues this iterative loop to refine the solution while monitoring convergence. Within the loop, several key things happen. (1) It applies (right) preconditioning, if necessary, (2) It computes a matrix-vector product to expand the Krylov space, and (3) It orthogonalizes the new vector against the recycle space as well as the Krylov space. Finally, (4) it computes the elements of the tridiagonal matrix, $T$, arising in the Lanczos iteration. After computing the QR factorization of $T$ to compute an update to the approximate solution $x_k$ at the $kth$ iteration, we must compute the eigenvalues of $T$, and update the solution. We have pointed out key steps in the algorithm that posed particular challenges when implementing the RMINRES algorithm.  There are other steps that are fairly straightforward and that we omit for brevity.  For instance, the code updates various variables for the next iteration, including values related to the Lanczos method. Upon convergence at iteration $m$, it returns the approximate solution $x_m$, as well as a number of other outputs, if requested by the user (e.g., a vector $relres$ containing the residuals at each iteration). It also provides the new recycle space to then use for the next linear system as the set of vectors we have referred to as $U$.

\subsection{Sparse Approximate Maps Implementation}\label{sec:SAMS}
We seek to eventually integrate SAMs as part of the larger predictor-corrector method code to explore their use as a mapping between matrices to recycle preconditioners when they are necessary. We again note that they were not required for the linear systems discussed in Figures \ref{fig:normDiff} - \ref{fig:MINRESvsRMINRES}, but in cases when the systems may become ill-conditioned, the role of SAMs will become vital to reduce costs associated with preconditioning many linear systems in a long sequence; see \cite{b7} for experiments demonstrating how SAMs reduced such costs. The update scheme in \cite{b7} solves 
\begin{equation}
N_k = arg min_{N\in \mathcal{S}} ||A_kN - A_0||_F
\label{eq:linear_system}
\end{equation}
and defines the updated preconditioner as $P_k = N_kP_0$.
Here, $\mathcal{S}$ represents the subspace defined by a selected sparsity pattern \cite{b7}. We assume the matrices $A_0$ and $A_k$ are part of a given sequence and are relatively close to each other.  

The code first allows the user to set a sparsity pattern, $\mathcal{S}$, that is imposed on the SAM. The code computes a map to transform one system matrix into another matrix in the sequence using equation (\ref{eq:linear_system}), which defines the full least squares problem (e.g., $n\times n$). However, with a judiciously chosen sparsity pattern (i.e., one that balances a small number of nonzeros with the accuracy as a map), the least squares equation given in (\ref{eq:linear_system}) results $n$ independent (and thus parallel), very small least squares problems. Specifically, it computes the least squares solution $n_j$ for the equation $Gn_j$ = $As_0$, where $G$ and $As_0$ are submatrices of $A_k$ and $A_0$, respectively, that respect the chosen sparsity pattern. Each solution, $n_j$, for the $n$ least squares problems defines the $j^{th}$ column of the approximate map $N_k$. 
The algorithms for SAMs and an analysis of the choice in sparsity pattern for several applications can be found in \cite{b7} (see Algorithms 4.1 and 4.2). 

\section{Challenges and Proposed Solutions}
We first discuss challenges we encountered while implementing the RMINRES algorithm in Python and PyTorch. Python requires a series of preprocessing steps to initialize appropriate data structures, variables, and parameters before beginning the algorithm directly.  One example is what happens after loading the files supplied (i.e., csv files storing the matrix elements), where the data must then be converted to arrays to perform calculations. Anther key difference is often the data must be reshaped to fit the variety of the SciPy sparse library and to suit user/application necessity (e.g., compressed sparse column or compressed sparse row).

\begin{figure}[h!]
    \centering
    \includegraphics[scale=0.45]{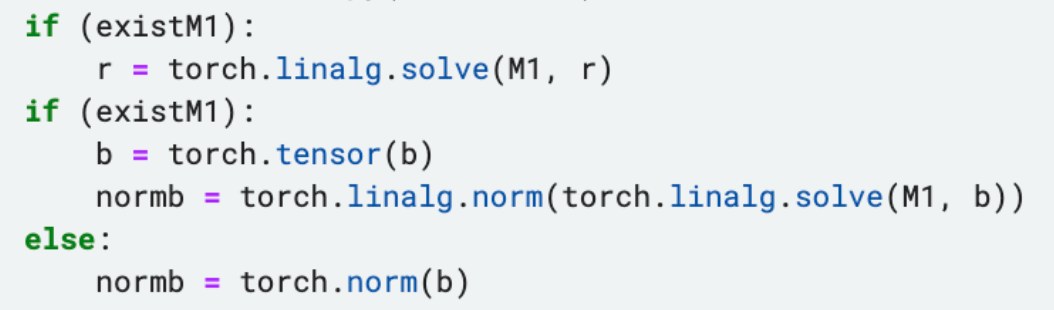}
    \caption{Code from RMINRES that fails due to underdevelopment of sparse tensors and torch.linalg.solve function.}
    \label{fig:A}
\end{figure}

We use compressed sparse column, and so the scipy.sparse function is used. The section of code seen in Figure \ref{fig:A} includes the first attempts at using torch.linalg.solve function for solving the systems of matrices. However, as seen in Figure \ref{fig:B}, the function does not execute. The error provided in Figure \ref{fig:B} is longer than shown below, however the most significant comments are provided. 

\begin{figure}[h!]
    \centering
    \includegraphics[scale=0.4]{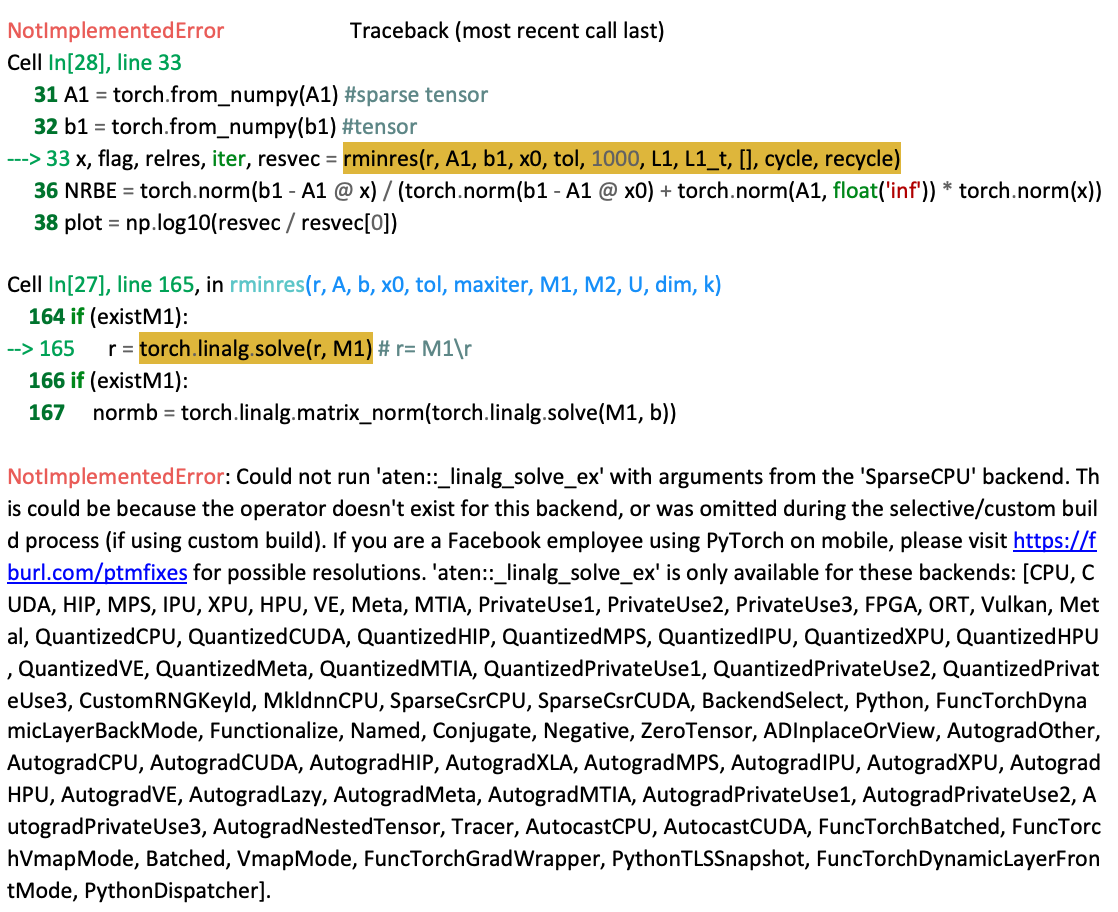}
    \caption{NotImplementedError after attempting to run torch.linalg.solve function in RMINRES implementation}
    \label{fig:B}
\end{figure}
The use of additional operations is fundamental to this predictor corrector method, as the matrices within this code must be in a format that allows for the solution of systems of linear equations and vector multiplication. In our implementations, it became clear that more development of the sparse tensor data structure is needed in order to successfully incorporate RMINRES into PyTorch.

Additionally, creating functions within PyTorch with diverse data type compatibility  will allow for solving a larger range of calculations within the PyTorch language and environment. This also includes different ways to compute matrix-vector multiplication for the predictor-corrector code, as $Hv$ is defined as a function and $H$ is not explicitly known. In the RMINRES code, we require not just matrix-vector multiplication, but matrix-matrix products (i.e, when orthogonalizing the subspace). More development is needed in this area, and part of ongoing work. While the predictor-corrector method and RMINRES solvers are implemented in Python \cite{b13}, the development of PyTorch sparse tensor data structure and compatible functions to complete the remaining portions of the RMINRES implementation is part of the ongoing work. 

\begin{figure}[h!]
    \centering
    \includegraphics[scale=0.5]{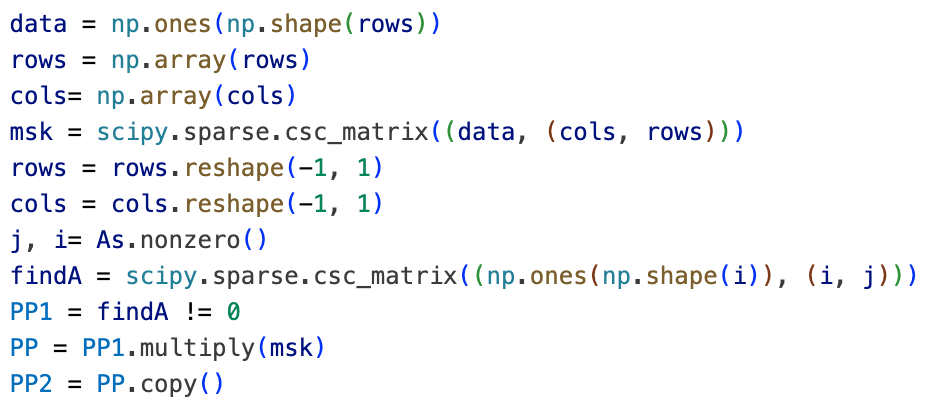}
    \caption{Examples of the data type agreement in Python.}
    \label{fig:C}
\end{figure}

Within the SAMs code, the issues we have encountered are related to mismatched data types, allocation of the appropriate data structure, and errors within called functions. As matrices are imported and converted to a sparse data structure, the updated matrices and data collected from computation needs to be converted and reshaped for the function at each step of the call. As with RMINRES, changes to the shape of the data are also required. Many of these issues have been resolved, however issues within the NumPy linalg.solve function called to solve linear systems are ongoing. An example of this is resolving data type matching errors for the implementation of the solver. Transitions between many data types occurs during the SAMs implementation that play a large role in the output. The segment of code shown in Figure \ref{fig:C} is part of the preprocessing of data for the functionality within the SAMs implementation. 

Lastly, to parallelize the code, we can utilize the Python multiprocessing library. This library allows us to create separate processes for individual tasks, which can run simultaneously and speed up the computation. We will divide the loop into separate tasks that can run concurrently. The code must be rewritten to reflect the parallel structure, particularly in the loop that computes the $n$ small least squares problems. Parallelizing the existing code is also part of ongoing work. 


\section{Discussion}

In comparing the computational costs expected from our ongoing implementation, an important note is assessing the efficiency gains offered by RMINRES in contrast to the conventional MINRES method. RMINRES has the potential to significantly reduce computational cost when solving sequences of linear systems by reducing the total number of iterations to converge to an accurate-enough approximate solution. This method updates iteratively to refine solutions, capitalizing on the shared data between consecutive linear systems. However, it is still important to note that the computational cost and time savings will depend on factors such as the size and sparsity of the linear systems, the nature of the coefficient matrices, and the specific problem at hand. SAMs can be useful when dealing with ill-conditioned matrices by reducing costs associated with computing a preconditioner from scratch for every linear system in the sequence. 



\section{Conclusions and Future Work}

We seek to reduce the overall cost of computing the linear systems arising in the predictor step of the predictor-corrector method utilized in \cite{b3} using recycling methods for both the iterative solver as well as for preconditioners. In this work, we considered the linear solver RMINRES and a map between coefficient matrices, the sparse approximate map. SAMs allow us to map closely related matrices, and may allow us to recycle previously computed preconditioners.  However, these methods do not yet exist in PyTorch, and thus have not yet been implemented in the Pareto front code described in \cite{b3}. While the predictor-corrector method and RMINRES solvers are implemented in Python, the development of PyTorch sparse tensor data structure and compatible functions is part of the ongoing and future work. With development in PyTorch linear algebra solver and sparse tensor data structure, we expect that the RMINRES and SAMs implementations will help to reduce computational cost in the overall predictor-corrector method for computing an optimal Pareto front\cite{b3}.

\section{Acknowledgments}

We thank the National Science Foundation for funding the Lehigh University Intelligent and Scalable Systems REU. Ainara Garcia is supported by the National Science Foundation under Grants NSF CNS-2051037.
Any opinions, findings, conclusions, or recommendations expressed in this document are those of the author(s) and should not be interpreted as the views of the National Science Foundation.

We also thank Sean Wang and Rishad Islam for their helpful consultations on this project.

\end{document}